\documentclass[letterpaper, 10 pt, conference]{ieeeconf}  

\IEEEoverridecommandlockouts                              
\usepackage{graphics} 
\usepackage{epsfig} 
\usepackage{amsmath} 
\usepackage{amssymb}  
\usepackage{subfigure}
\usepackage{color}
\title{\LARGE \bf
Statistical tests for group comparison of manifold-valued data}

\author{Anne Collard, Christophe Phillips and Rodolphe Sepulchre
\thanks{This paper presents research results of the Belgian Network DYSCO (Dynamical Systems, Control, and Optimization), funded by the Interuniversity Attraction Poles Programme, initiated by the Belgian State, Science Policy Office. The scientific responsibility rests with its author(s).}
\thanks{A. Collard is with the Departement of Electrical Engineering and Computer Science,
    University of Li\`{e}ge, 4000 Li\`{e}ge, Belgium;       
     {\tt{anne.collard}@ulg.ac.be}}%
\thanks{R. Sepulchre is with he Departement of Electrical Engineering and Computer Science,
    University of Li\`{e}ge, 4000 Li\`{e}ge, Belgium, and INRIA Lille-Nord Europe, Orchestron project.
40 avenue Halley
F 59650 Villeneuve d'Ascq, France;    
     {\tt{r.sepulchre}@ulg.ac.be}}%
\thanks{C. Phillips is with the Cyclotron Research Centre, 
    University of Li\`{e}ge, 4000 Li\`{e}ge, Belgium;
        {\tt{c.phillips}@ulg.ac.be}}%
}

\begin{document}

\maketitle
\thispagestyle{empty}
\pagestyle{empty}

\begin{abstract}
Motivated by population studies of Diffusion
Tensor Imaging, the paper investigates the use of
mean-based and dispersion-based permutation tests
to define and compute the significance of a
statistical test for data taking values on nonlinear manifolds.
The paper proposes statistical tests that are computationally
tractable and geometrically sound for Diffusion Tensor Imaging.
\end{abstract}

\section{INTRODUCTION}
Statistical analysis of scalar-valued  images is a well established central component
of contemporary science.  But the evolution of sensor technology and data storage
increasingly produces images that  are multimodal and nonlinear in nature. This 
has motivated significant work in the recent years to extend signal processing techniques 
from scalar-valued data to manifold-valued data, see e.g. \cite{jensen1996introductory}.

The present paper is specifically motivated by the   increasing role  of diffusion tensor imaging 
in neuroscience. In their
simplest form, DT images provide a $3 \times 3$ positive definite diffusion tensor for each voxel \cite{Basser1994259}.
Voxel-based statistical analysis of a population therefore involves statistical tests
among positive definite tensors rather than scalars. The default remedy is to convert the
tensor information into a scalar information (usually fractional anisotropy, see below), but
both the intensity information (that is, the three positive eigenvalues of the tensor) and
the orientation information (the orientation of the three principal axes) potentially contain
valuable statistical information, calling for new methodological developments.

The challenge is methodological as well
as computational because clinical studies usually involve large populations and many
voxels, that is, large-scale statistical analysis.

With this motivation in mind, the present paper investigates the methodological and 
computational value of two standard non-parametric permutation tests: a mean-based permutation test
and a dispersion-based permutation test. We discuss how these tests can be extended from
scalar-valued data to manifold-valued data and specialize the discussion in the particular
case of Diffusion Tensors, that is, $3 \times 3$ positive definite matrices.

We argue that a dispersion-based permutation test is a computationally tractable approach
to clinical DTI statistical analysis and stress the value of defining properly motivated geometric
quantities for the underlying similarity measures. 

\subsection{State of the art for statistical analysis of DTI}
There is a vast recent literature on DTI statistical analysis. Among these papers, we can distinguish the one using univariate tests to compare DT images and the ones trying to exploit further information. For example, \cite{Chappell08}, \cite{Chung08}, and \cite{Jones05} focused on the fractional anisotropy or the mean diffusivity of tensors, without paying any attention to the orientation of these tensors. These scalars are indeed invariants linked to the shape of tensors, but can not detect any difference in orientation between tensors. In \cite{Goodlett09}, the evolution of the fractional anisotropy along a fiber tract is studied. Some other papers have tried to use the whole information contained in the tensor, namely through the use of the Log-Euclidean metric (\emph{i.e.} they use the logarithms of tensors instead of the tensors themselves). This is the case in \cite{Lee07}, where an Hotelling's $T^2$ test is developped for the Log-Euclidean metric, a framework similar to the ones in \cite{Fouque11} and in \cite{Grigis12}. Some other papers develop a rigorous conceptual framework based on the Riemannian manifold $\mathrm{S}^+(3)$, as in \cite{Pennec06}, \cite{Fletcher07riemanniangeometry} and \cite{Lenglet:2006uq}. However, those papers have not addressed the statistical significance of a test for group comparison. This is also the case of \cite{Dryden09}, which uses another parametrization of tensors. Statistical tests were proposed in \cite{Schwartzman10}, through a decoupled analysis of the eigenvalues and the eigenvectors of the tensors. The tests are based on distributional assumptions, which is a potential limitation for diffusion tensors. The closest published work to the present paper is \cite{Whitcher07}, when the authors propose a multivariate dispersion-based permutation test, see Section III for more details.\\

The paper is organized as follows: after a brief review of the state of the art for the statistical analysis of Diffusion Tensor, Section II will focus on statistical tests for groups comparison, beginning with mean-based permutation tests. Computation of means on manifolds will also be discussed, before the introduction of dispersion-based tests. The case of multivariate tests will be addressed, and some appropriate similarity measures for DTI will be introduced. Section III deals with the methods that we have used for our tests, while Section IV shows several of our results.

\section{STATISTICAL TESTS FOR GROUPS COMPARISON}
Statistical analyses of scalar images are often performed through the use of parametric methods \cite{Friston07}, such as the Student $t^2$ test for comparison of Gaussian variables. However, the distribution of multivariate data is rarely known and, if known, is often not Gaussian. This explains why many authors have made the choice of non-parametric methods to study multivariate images. Among these non parametric methods, permutation tests are often used because of their relative simplicity. Permutation methods provide statistical significance testing of difference between groups without having to assume a distribution of the data. These methods have the ability to directly estimate the null distribution of the statistics describing the difference. Moreover, these methods are easily applicable to any statistical test, which is interesting to compare results obtained with different parametrizations of the data. 
\subsection{Mean-based permutation tests}
\begin{figure}[thpb]
      \centering
      \includegraphics[scale=0.31]{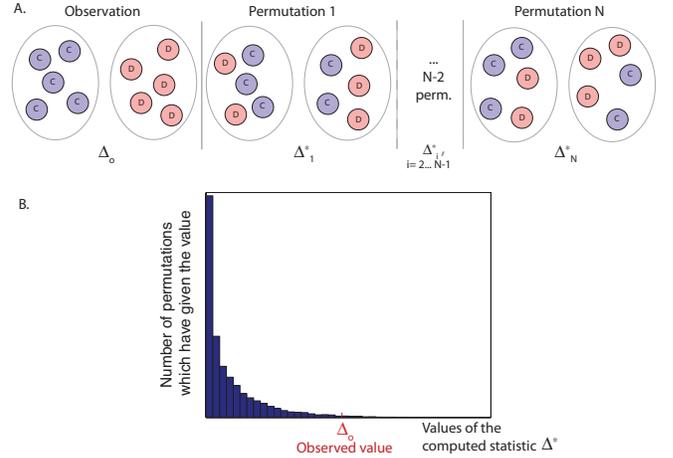}
      \caption{General procedure of a permutation test. Given the observed value of the statistic, values corresponding to $N$ permutations of the labels are computed. This gives an approximation of the (unknown) distribution of the statistic. The comparison between the observed value and this distribution enables to compute the p-value of the test. }
      \label{perm_expl}
   \end{figure}

Permutation tests are based on a simple idea. For the sake of illustration, consider the statistical significance of a variable $x$
to distinguish among two populations C and D. Suppose that the mean $\mathbb{E}(x)$
differs by a quantity $\Delta_0$ between group C and group D. Permutation
tests enable to quantify, without assumptions about the distribution of the variable, if this difference is significant or not. Indeed, if the difference is not significant,
it should not be altered by random permutations between C and D. Therefore,
given the null hypothesis that the labelings are arbitrary, the significance of
$x$ can be assessed by comparison with the distribution of values under all possible
permutations. This is illustrated by the histogram in Figure \ref{perm_expl}. If the observed
difference $\Delta_0$ is in the tail of the distribution, it means that very few permutations
of the data attain the same difference. The $p$-value of the test is given by the ratio between the number of times that a permuted statistics is higher than the observed value and the number of performed permutations. The test is statistically significant at a level $\alpha$ if the $p$-value is smaller than $\alpha$. This means that this value has less than $100 \alpha\%$ of chance to have be found randomly. 

The generalization of a permutation test to data that takes values on a manifold is
conceptually straightforward because it only requires a proper notion of mean.
On a Riemannian manifold $\mathcal{M}$, the Karcher mean of a set of points $\{x_1, x_2, \dots x_N\}$ is given 
by the Fr\'echet formula
\begin{equation} \mu = \mathrm{arg} \min_{x \in \mathcal{M}} \frac{1}{2N} \sum_{i=1}^N d^2(x, x_i) \, , \end{equation} (with $d$ the distance on the manifold), which reduces to the classical arithmetic mean when
using the Euclidean distance.

\subsection{Means on manifolds and means for DTI}
Computing Riemannian means (or medians) on specific manifolds has been the object of
significant research in the recent years, see e.g. \cite{Absil04} for means on the Grassmann manifold,
\cite{Pennec06b,Fletcher07riemanniangeometry} for means and medians on the space of positive definite matrices, and \cite{Bonnabel:2009uq} for a mean
for fixed-rank positive semidefinite matrices.

It should be emphasized that a mean-based permutation test may represent a formidable computational
task. For instance, computing the Riemannian mean of several positive definite matrices is typically
achieved by an iterative algorithm \cite{Pennec06b}.  For a population of size $N$ classified in $g$ groups of size $n_i$, this computation must be repeated 
$ M= N! / \prod_{i=1}^g n_i!$  which is prohibitive for large populations, even if the full distribution is not computed.
In the context of DTI, the statistics are typically computed for a large number of voxels, which adds
to the computational burden.

One computational remedy for positive definite tensors is to compute the arithmetic mean of tensors
(that is, to work with the Euclidean metric) or, better,  to compute a matrix  geometric mean according
to the formula
\begin{equation} \mu_{\mathrm{LE}} = \exp(\sum_{i=1}^N \log x_i) \, .\end{equation}
which corresponds to the Riemannian mean for the Log-Euclidean metric \cite{Arsigny:2006fk}.

In a recent paper, the authors have introduced another notion of mean for positive tensors,
that provides a better decoupling between orientation and anisotropy \cite{Collard12}. Using the spectral
decomposition  of each tensor, the mean tensor is defined as 
\begin{equation} \mu_{SQ} = R_{\mu}\Lambda_{\mu} R_{\mu}^T\, , \end{equation} 
where $\Lambda_{\mu}$ is the geometric mean of the eigenvalues of the tensors, \emph{i.e.} $\Lambda_{\mu} = \exp(\sum_{i=1}^N \log \Lambda_i)$, where $\Lambda_i = (\lambda_{1,i},\lambda_{2,i}, \lambda_{3,i})$, with the ordered eigenvalues $\lambda_{1,i}>\lambda_{2,i}> \lambda_{3,i}$. The mean orientation $R_{\mu}$ is computed through the chordal mean of quaternions.

Even with such computational simplifications, the computational burden of a mean-based
approach remains prohibitive for a real application of group studies using Diffusion Tensor
Images because the matrix $\log$ operation will scale in a factorial way with  the population size. 

\subsection{Dispersion-based permutation tests}
The Multiresponse permutation procedure (MRPP) proposed in \cite{Mielke01} is not based on
repeated computations of means but only requires to build a similarity matrix between
all data points. Given a symmetric similarity measure $s(i,j)$ between data points $x_i$
and $x_j$, the similarity matrix is defined as the symmetric matrix $S$ with element $S_{ij}= s(i,j)$.

The statistical test proposed in \cite{Mielke01} is based on a measure of dispersion of the data points
within each group rather than on a measure of mean. The dispersion $\delta_i$ of a group
of $n_i$ elements is defined as
\begin{equation} \delta_i = \frac{2 (n_i-2)!}{n_i!} \sum_{i<j} s(i,j)\end{equation}
where the sum is computed over all data points of the group.

The overal dispersion $\delta$ of the variable $x$ in the population is defined as a weighted
sum of the dispersions in each group:
\[ \delta= \sum_{i=1}^g C_i \delta_i \]
where $C_i>0, i = 1 \dots g$  are the weights of each group (their sum must be 1).

The rest of the procedure follows the permutation test described in Section II. A. , with the mean-difference
$\Delta$ replaced by the dispersion $\delta$. The observed dispersion is judged statistically significant
only if it occurs in the lower tail of the histogram among all possible permutations of the population.

The distance-based permutation test has the same advantages than the mean-based permutation
test: it does not require any assumption about the statistical distribution and only requires
a similarity measure between data points, for instance a distance on Riemmannian manifolds.
But a significant computational advantage is that the similarity matrix must be computed only
once, requiring $O(N^2)$ computations of pairwise similarity for a total population of size $N$.

\subsection{Multivariate testing}

Our description of permutation tests has so far assumed an univariate statistical testing but
is easily extended to multivariate testing. A widely used method for the statistical comparison of multivariate data consists in the computation of 'marginal` or 'partial` tests (one for each of the $n_v$ considered variable) and to combine the $p$-values obtained for each partial test by a combining function, which can be of different forms \cite{Pesarin01}. This method involves two stages of computation. First, the marginal $p$-values (denoted $\xi_i, i = 1, \dots, n_v$) are computed through a permutation test (here, we will use a MRPP test for univariate data). Then, the combining function is used to compute the combined observed value, $T_o = \mathcal{C}(\xi_1, \dots, \xi_{n_v})$. The distribution of this $T$ is computed through the combination of the $p$-values computed for each permutation of the first step, \emph{i.e.} $T^*_k = \mathcal{C}(\xi^*_{1,k}, \dots, \xi^*_{n_v,k})$. The combined $p$-value for the test is then estimated via the ratio between the number occurrences where $T^*_k$ was larger than $T_o$ and the total number of performed permutations.

\subsection{Similarity measures for DTI}
Every distance on the manifold of positive tensors qualifies for a similarity measure. In the
present work, we will compare two similarity measures to the naive Euclidean distance between matrices: the Log-Euclidean distance \cite{Arsigny:2006fk} and the spectral-quaternion similarity measure recently introduced in \cite{Collard12}. 
Both similarity measures only involve a spectral decomposition of a $3\times 3$ matrix
for all  data points as the main computation of the similarity measure.  

Regarding multivariate testing, we propose to compare an "Euclidean" statistical test
based on the six independent quantities of a $3 \times 3$ positive definite tensor (as proposed
in \cite{Whitcher07}) and a "geometric" statistical test based on the six geometric quantities that define
scaling and orientation of the tensor: the three eigenvalues and the three
first components of the quaternion. In the latter case, we use a geometric
similarity measure for the (positive) eigenvalues $s(\lambda, \mu) = \sqrt{\log^2\left( \frac{\lambda}{\mu} \right)}$
and an Euclidean (chordal) measure for the quaternions. A further alternative would be a 'log-euclidean` statistical test based on the six elements of the logarithm of the tensor, as investigated in \cite{Whitcher07}. We do not include this comparison here since the work of \cite{Whitcher07} suggests that it does not offer significant advantages compared to the Euclidean test.

For comparison purposes, we will also compute an univariate test using the fractional
anisotropy of tensors \cite{Chung08}.
\section{METHODS}
We will compare the power of the proposed statistical tests on the following synthetic data sets.

We generate two groups of tensors starting from a reference tensor  and a transformed tensor. 
The parameter $\gamma$ quantifies the amount of deformation. Denoting by $\lambda_i, i = \{1,2,3\}$, the
eigenvalues of the reference tensor and by $\lambda'_j$ the ones of the deformed tensor, the following four different geometric transformations
of the reference diffusion tensor are studied \cite{Boisgontier12}:
\begin{itemize}
\item \emph{Decrease of longitudinal diffusion (DL)}: $\lambda'_1$ will be given by $\lambda_1' = \lambda_1- \gamma(\lambda_1-\lambda_2),\, \gamma \in [0,1]$. The two others eigenvalues will be left unchanged.
\item \emph{Increase of radial diffusion (IR)}: $\lambda_2$ and $\lambda_3$ will be replaced by $\lambda'_j = \lambda_j + \gamma \frac{\lambda_1-\lambda_j}{2},\, j= \{2,3\},\, \gamma \in [0,1]$, while $\lambda_1$ will be left unchanged. 
\item \emph{Increase of mean diffusion (IM)}: all the eigenvalues of the deformed tensor will be given by $\lambda'_j = (1+\gamma) \lambda_j, \, j = \{1,2,3\}, \, \gamma \in [0,1]$.
\item \emph{Change of diffusion orientation (CO)}: the angle $\theta$ between principal directions of the reference tensors will change following $\theta' = \theta + \gamma \frac{\pi}{2}, \, \gamma \in [0,1]$.
\item \emph{Change in both eigenvalues and orientation}: in this case, we will combine a difference in orientation (CO) with one of the first three differences. Both modifications will evolve following the same $\gamma$.
\end{itemize}

Starting from the reference tensor and its reference deformation, we generate a population
of $N = 20$ tensors by sampling from a Wishart distribution 10 tensors  around the reference tensor
and 10 tensors around the deformed tensor. 

The statistical comparison will be tested in a situation of high anisotropy, where the eigenvalues of the reference are $\Lambda = (5,1,0.5)$, a situation of small 
anisotropy where they are given by $\Lambda = (3, 1, 1)$, and a situation of near isotropy, $\Lambda= (1.3, 1, 1)$.
The performance of the different tests in many situations will be assessed through the computation of the \emph{power} of these tests. This quantity is computed by performing the same test a large amount of times and by counting the occurrences of significant results. The power is the ratio between this number of occurrences and the total amount of performed tests. If the test is not efficient, the power will be close to the significance level $\alpha$. A value of the power equal to one means that the test is particularly efficient for the analyzed situation. The different used values (number of tests, level of significance, \dots) are summarized in Table \ref{params}. It should be noted that the level of noise is relatively high for these tests. The tests could be done with a larger number of degree of freedom in the Wishart distribution (which corresponds to a lower level of noise).

\begin{table}
\caption{Parameters of the statistical tests}
\label{params}
\begin{center}
\begin{tabular}{c|c}

Parameter & Value\\
\hline
Level of significance $\alpha$ & 0.05 \\
Number of samples by groups $n_i$ & 10 \\
Number of permutations $N_p$ & 20 000\\
Number of tests by situation $N_t$ & 500\\
\end{tabular}
\end{center}
\end{table}

\section{RESULTS}
In the following, some interesting results found during our tests will be shown. We will first study how dispersion-based tests vary with the used similarity measure and then observe the difference in interpretation of the results which can be done from multivariate tests. Due to space constraints, other results will not be shown.

\subsection{Univariate tests}
\begin{figure}[thpb]
      \centering
      \includegraphics[scale=0.45]{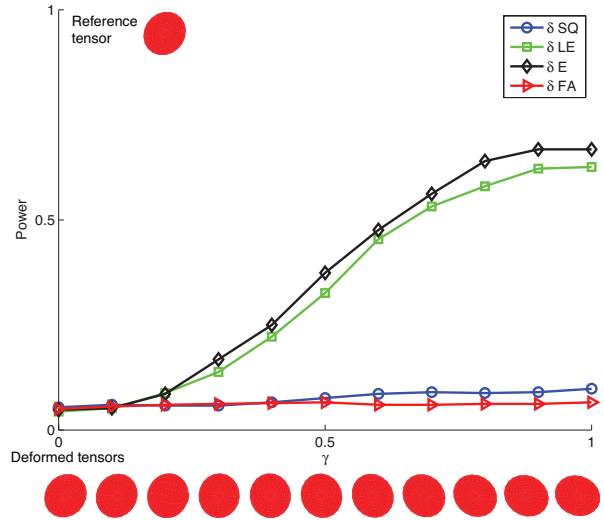}
      
      \caption{Power of the statistical tests in a situation of very low anisotropy (the tensors are nearly isotropic). The tests use the spectral quaternion measure ($\delta$ SQ, blue circles), the Log-Euclidean measure ($\delta$ LE, green squares), the Euclidean measure ($\delta$ E, black diamonds) and the Fractional Anisotropy of tensors ($\delta$ FA, red triangles). A difference in the diffusion orientation (CO) is simulated. In this case, as illustrated on the figure, the tensors are very similar (all almost spheric) and we argue here that no difference should be noted. However, this is not the case for the Log-Euclidean and Euclidean tests. }
      \label{p_tests_CDC2c}
   \end{figure}

Figure \ref{p_tests_CDC2c} illustrates an uncommon situation, where the desired result is not a maximal power of the test, but a minimal one. The simulated situation is the one of near isotropic tensors, with a progressive change of orientation. As the tensors are all nearly spheric, there should be no noticed difference between them, and the tests should not be sensitive to the simulated deformations. However, it can be seen that this is not the case using Euclidean or Log-Euclidean measures. In this case, the Spectral Quaternions measure and the test based on Fractional Anisotropy perform better, as they do not detect any difference between the reference and the deformed tensors. It is interesting to note that these curves of power can also be interpreted as curves of sensibility and robustness of the measures. The power of a test is a measure of its sensitivity to the considered deformation. Conversely, flat curves indicate a robustness of the test to a given deformation.
    \begin{figure}[thpb]
      \centering
      \includegraphics[scale=0.45]{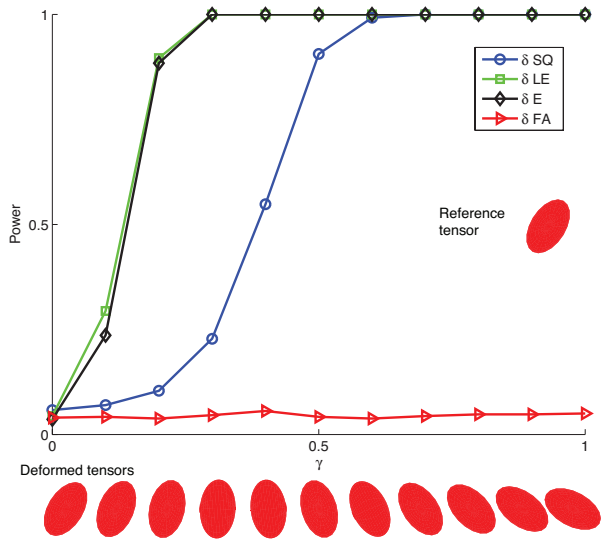}
      
      \caption{Power of the statistical tests in a situation of low anisotropy. A difference in the orientation of the tensors is simulated. The tests use the spectral quaternion measure ($\delta$ SQ, blue circles), the Log-Euclidean measure ($\delta$ LE, green squares), the Euclidean measure ($\delta$ E, black diamonds) and the Fractional Anisotropy of tensors ($\delta$ FA, red triangles). For $\gamma = 0$, which means no difference between the two references of the groups, the power of the tests is about $\alpha$. As expected, the test based on the Fractional Anisotropy fails to detect the difference in orientation. The Euclidean and Log-Euclidean tests are very sensitive to this kind of deformation, while the Spectral Quaternion measure is between those two situations.}
      \label{p_tests_CDC2b}
   \end{figure}

This analysis of robustness is also applicable to the example of Figure \ref{p_tests_CDC2b}, showing a difference in orientation in the case of low anisotropy. In this figure, a clear difference can be observed between the Fractional Anisotropy test and the Euclidean and Log-Euclidean tests. As expected, the Fractional Anisotropy test is never significant for this type of difference, showing one more time the limitations of using a unique scalar to represent multivariate data, as the orientation information is totally lost. To the contrary, the Euclidean and Log-Euclidean tests are very sensitive to this difference and thus exhibit strong performance. The Spectral Quaternion measure offers an intermediate situation: it is sensitive to the orientation information, but less that the two Euclidean tests. In fact, the orientation term of this measure is weighted by a parameter $k$, which depends on the anisotropies of tensors. The role of this parameter is to decrease the importance of the orientation term in case of low anisotropy, since in this case, the orientation information is highly uncertain. This explains why the Spectral Quaternion measure is not sensitive to the deformations shown in Figure \ref{p_tests_CDC2c}. It is important to understand the impact of this parameter of the results. If it is increased, the orientation term will become more important, which will produce an increase of the sensitivity of the measure (\emph{i.e.} an increase of the performance of the test). Depending upon the tradeoff between sensitivity and robustness, the curve of the Spectral Quaternion test can be closer to the Fractional Anisotropy one, or to the contrary, closer to the Log-Euclidean test. The fact that this tradeoff can be tuned is of relevance for clinical applications.
   
\subsection{Multivariate tests}   
\begin{figure}[thpb]
      \centering
      \subfigure[Geometric parametrization]{\includegraphics[scale=0.42]{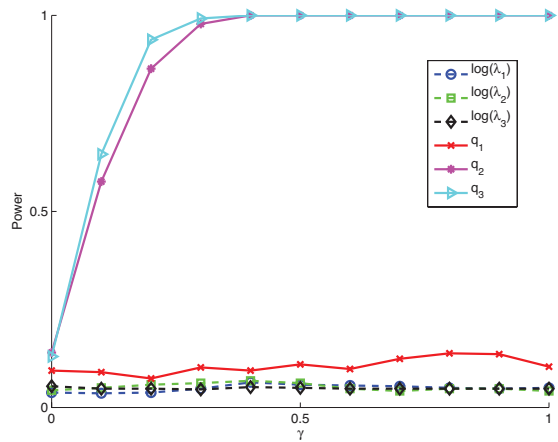}}
       \subfigure[Euclidean parametrization]{\includegraphics[scale=0.42]{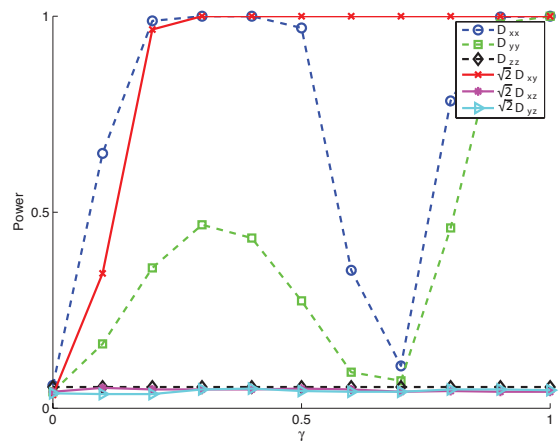}}
      
      \caption{Power of each of the partial statistical tests for multivariate parametrizations of the tensors. (a) A geometric parametrization is used. (b) An Euclidean parametrization is used. The simulated situation is a deformation of orientation in a case of high anisotropy. The results of the geometric tests are easily interpretable, as only the partial tests associated to the orientation are performant. To the contrary, the Euclidean tests are poorly understandable.}
      \label{p_tests_CDC_multi2}
   \end{figure}

\begin{figure}[thpb]
      \centering
      \subfigure[Geometric parametrization]{\includegraphics[scale=0.42]{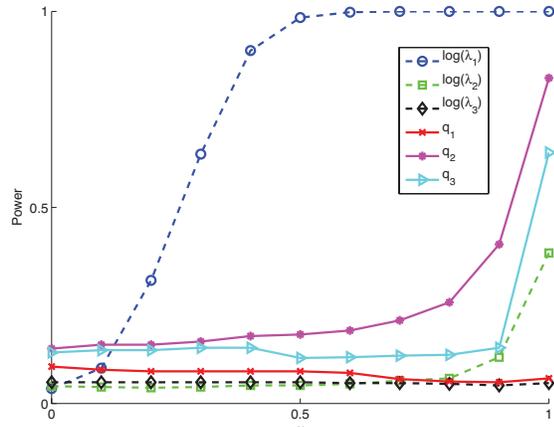}}
       \subfigure[Euclidean parametrization]{\includegraphics[scale=0.42]{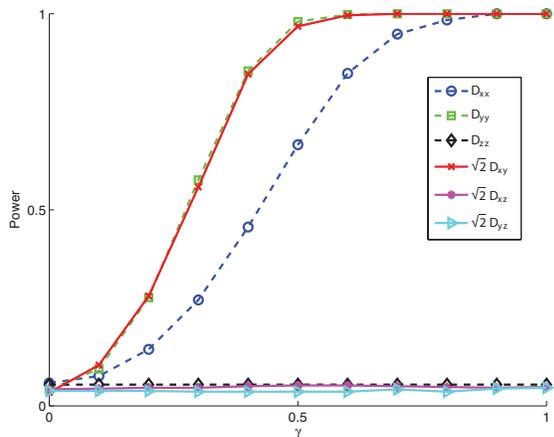}}
      
      \caption{Power of each of the partial statistical tests for multivariate parametrizations of the tensors. (a) A geometric parametrization is used. (b) An Euclidean parametrization is used. The simulated situation is a decrease of the longitudinal diffusion in a case of high anisotropy. The results of the geometric tests are easily interpretable, as only the partial test associated to the first eigenvalue is performant. To the contrary, the Euclidean tests are difficult to interpret.}
      \label{p_tests_CDC_multi1b}

   \end{figure}

In the following, we will focus on the interpretation of the results of multivariate tests. \\
Figure \ref{p_tests_CDC_multi2} illustrates the results of each partial tests for multivariate parametrizations of the tensors, for a simulated change of orientation in a case of high anisotropy. The comparison between geometric parametrization (\emph{top}) and algebraic one (\emph{bottom}) is straightforward. As the geometric parametrization clearly shows that the orientation only has been changed, this interpretation can not be drawn from the results of the Euclidean tests. This information could however be of great importance in clinical studies. In a similar way, the decrease of longitudinal diffusion simulated in Figure \ref{p_tests_CDC_multi1b} is clearly seen with the geometric parametrization, while this is not the case for the Euclidean one. Indeed, for a geometric parametrization, the partial test of the first eigenvalue is the only one to detect a difference. It should be noted that, from $\gamma = 0.9$, the first eigenvalue is very close to the second one, which increases the uncertainty in orientation (this explains why other partial tests become significant). The easy interpretation of the statistical tests using the geometric parametrization is a desirable feature, which opens the way to many applications. It should be noted that if a unique significance level is needed, the partial tests can be combined using an appropriate function.

%

\section{CONCLUSIONS}

In this work, we have shown that existing methods in the field of groups comparison could be advantageously used for the statistical analyses of data lying on Riemannian manifolds. These methods have several advantages as they do not use any assumptions about the distribution of the data (which is seldom known). Moreover, it has been shown that the only specific tool which is needed for this group comparison is an appropriate similarity measure between the data (or a parametrization of them). We have illustrated the computational advantage of basing the permutation test on dispersion rather than means.\\
Using Diffusion Tensor Images as an example, we have shown how different measures or different parametrizations of the data can affect the results of the tests. Moreover, two interesting features of the spectral quaternion measure (and the geometric parametrization associated to this measure) have been highlighted, which both could be relevant for clinical applications.\\
Both for computational and conceptual reasons, dispersion based permutation tests offer an appealing framework for group comparison of manifold-valued data.


\bibliographystyle{IEEEtran}      
\bibliography{biblio_cdc}   






\end{document}